\documentclass[12pt,leqno]{article}
\usepackage{amsmath}
\usepackage{amssymb}
\usepackage{theorem}
\usepackage{pstricks,pst-plot}
\usepackage{mathrsfs}

\font\teneufm=eufm10 scaled \magstep1
   \font\seveneufm=eufm7 scaled \magstep1
   \font\fiveeufm=eufm5 scaled \magstep1
  \newfam\eufmfam
  \textfont\eufmfam=\teneufm 
   \scriptfont\eufmfam=\seveneufm
   \scriptscriptfont\eufmfam=\fiveeufm

\def\hexnumber@#1{\ifcase#1 0\or 1\or 2\or 3\or 4\or 5\or 6\or 7\or 8\or
 9\or A\or B\or C\or D\or E\or F\fi}

\newcounter{newsection}[section]
\theoremstyle{change}
\newtheorem{thm}{Theorem.}[newsection]

\newtheorem{lem}[thm]{Lemma.}

\newtheorem{cor}[thm]{Corollary.}
\newtheorem{example}[thm]{Example.}

{\theorembodyfont{\normalfont\rmfamily}
}

{\theorembodyfont{\normalfont\rmfamily}
}

\makeatletter

\newcommand{\proofname}{Proof}

\renewcommand\section{\@startsection{section}{4}{\z@}%
                                    {3.25ex \@plus1ex \@minus.2ex}%
                                    {-1em}%
                                    {\reset@font\normalsize\bfseries}}
\renewcommand\subsection{\@startsection{subsection}{5}{\z@}%
				    {0ex \@plus1ex \@minus.2ex}%
                                    {-1em}%
                                    {\reset@font\normalsize\bfseries}}
\makeatother

\input cyracc.def

\font\tenscr=rsfs10  scaled \magstep1

\font\sevenscr=rsfs7  scaled \magstep1
\font\fivescr=rsfs5  scaled \magstep1
\skewchar\tenscr='177 \skewchar\sevenscr='177 \skewchar\fivescr='177
\newfam\scrfam \textfont\scrfam=\tenscr \scriptfont\scrfam=\sevenscr
\scriptscriptfont\scrfam=\fivescr
\def\scr{\fam\scrfam}

\numberwithin{equation}{section}

\def\sC{{\scr C}}

\def\sO{{\scr O}}\def\sP{{\scr P}}
\def\sR{{\scr R}}

\def\bC{\mathbb C}

\def\bR{\mathbb R}\def\bS{\mathbb S}\def\bT{\mathbb T}
\def\bU{\mathbb U}
\def\bZ{\mathbb Z}

\def\ve{\varepsilon}
\def\vp{\varphi}

\def\vt{\vartheta}

\def\what{\widehat}

\def\vep{\varepsilon}

\def\z{\zeta}

\def\cn{\bC^n}

\def\bvt{ {\rm{CBV}}^n(\bT)}
\def\gammasq{\gamma\times\gamma}

\def\bvj{ {\rm{CBV}}^n(J)}
\def\bvjr{ {\rm{CBV}}^n_{\bR}(J)}

\def\bv{\mathbf v}

\def\Z{\mathbb Z}

\def\T{\mathbb T}

\def \R{\mathbb R}

\def\wg{\widehat\gamma}

\def\c2{\mathbb C^2}

\def\sqgamma{\{\gamma_k\}_{k=1,2,\dots}}

\def\sqhatgamma{\{\what\gamma_k\}_{k=1,2,\dots}}

\def\sqsigma{\{\sigma_k\}_{k=1,2,\dots}}

\def\Rm{\begingroup\let\par=\null\obeylines\RemovE}
\def\RemovE#1\mR{\endgroup}

\begin{document}
\ 

\vskip 1in

\centerline{\bf T{\sc{he}} C{\sc{onvergence of}} H{\sc{ulls of}} C{\sc{urves}}}
\bigskip

\centerline{Alexander J. Izzo$^($\footnote{This author was partially supported by a Simons collaboration grant and by NSF Grant DMS-1856010.}$^)$ and Edgar Lee Stout}

\begin{abstract}
It is shown that a simple closed curve in $\cn$ that is a uniform limit of rectifiable simple closed curves each of which has nontrivial polynomial hull has itself nontrivial polynomial hull.  In case the limit curve is rectifiable, the hull of the limit is shown to be the limit of the hulls.  It is also shown that every rectifiable simple closed curve in 
$\cn$, $n\geq 2$, can be approximated in total variation norm by a polynomially convex, rectifiable simple closed curve that coincides with the original curve except on an arbitrarily small segment.  As a corollary, it is shown that every rectifiable arc in $\cn$, $n\geq 2$, is contained in a polynomially convex, rectifiable simple closed curve.
\end{abstract}

\vskip .3in

%
%

\section{Introduction.} 
It is known from the work of Forstneri\v{c}  \cite{Forstneric:1994, Forstneric:1986}, Forstneri\v{c} and Rosay \cite{Forstneric-Rosay:1993}, and L{\o}w and Wold \cite{Low-Wold:2009} that given a compact, smooth manifold $M$ of dimension $d<n$, the set of polynomially convex, totally real embeddings of $M$ in $\cn$ of class $\sC^s$, $1\leq s\leq\infty$, is open and dense in the space of all embeddings of $M$ in $\cn$ of class $\sC^s$ in the $\sC^s$ topology.  Specialized to simple closed curves, this says that the set of polynomially convex simple closed curves of class $\sC^s$ in $\cn$, $n\geq2$, is open and dense in the space of all simple closed curves in $\cn$ of class $\sC^s$ with the $\sC^s$ topology. 

In this paper we establish stronger results regarding openness and denseness of polynomially convex simple closed curves that, in particular, contain the statement that the polynomially convex, rectifiable simple closed curves in $\cn$, $n\geq2$, form a dense open set in the space of all rectifiable simple closed curves with respect to both the supremum norm and the total variation norm.  
Regarding openness, we show that if a simple closed curve $\gamma$ is a uniform limit of rectifiable simple closed curves each of which has nontrivial polynomial hull, then $\gamma$ has nontrivial polynomial hull.  No regularity hypothesis is made on the limit curve $\gamma$ in this result.  We also show that when the limit curve $\gamma$ is rectifiable, the hull of the limit is the limit of the hulls, but this can fail for non-rectifiable $\gamma$.  Regarding denseness, we show that every rectifiable simple closed curve in $\cn$, $n\geq 2$, can be approximated in total variation norm by a polynomially convex, rectifiable simple closed curve that coincides with the original curve except on an arbitrarily small segment.  Analogous results concerning approximation of smooth simple closed curves in $\sC^s$ topologies are also presented.

As a corollary of our density results we prove that every rectifiable arc in $\cn$, $n\geq 2$, is contained in a polynomially convex, rectifiable simple closed curve, and every $\sC^s$-smooth arc is contained in a polynomially convex, $\sC^s$-smooth simple closed curve.  It is also true that \emph{every} polynomially convex arc is contained in a polynomially convex simple closed curve.  However, the proof of that result will be published separately.

In the next section we introduce some terminology and notation.  Our main results are stated in Section~\ref{results}  Proofs and related lemmas and examples are given in Sections~\ref{convergence}--\ref{density}

%
%

\section{{\bf {Terminology and Notation}}.}\label{notation}
\setcounter{newsection}{2}
\setcounter{thm}{0}

We denote the topological boundary of a subset $B$ of $\cn$ by $bB$.
Neighborhoods will always be taken to be open sets.
We use the standard notation $\|f\|_{X} = \sup\{ |f(x)| : x \in X \}$ for the supremum of a bounded complex-valued function $f$ over a set $X$.  If $X$ is a compact subset of $\cn$, the \emph{polynomial hull} 
$\what X$ of $X$ is defined by
$$\what X=\{z\in\cn:|P(z)|\leq \| P \|_X\ 
\mbox{\rm{for\ all\ polynomials}}\ P\}.$$
The set $X$ is said to be \emph{polynomially convex} if $\what X = X$.  The polynomial hull of $X$ is said to be \emph{nontrivial} if instead the set $\what X \setminus X$ is nonempty.

Let $J$ be either a closed interval in the real numbers or a circle.  We 
denote by $\bvj$ the space of all continuous maps from $J$ into $\cn$ that are of bounded variation.  The corresponding space of maps from $J$ into $\bR^n$ will be denoted by $\bvjr$.  We denote the total variation of a map $F$ over $J$ by ${\rm{var}}\,F$.
The spaces $\bvj$ and $\bvjr$ are Banach spaces with the total variation norm $\|\cdot\|_{\rm bv}$ given by
 $$\|F\|_{\rm bv}=\| F \|_J+{\rm{var}}\, F.$$
In case $F$ is injective, the total variation of $F$ over $J$ is simply the length of the image $F(J)$.
We denote the unit circle by $\bT$.
The elements of $\bvt$ are then the \emph{rectifiable closed curves}. A closed curve is {\it{simple}} if it is injective. Frequently below we will use the common abuse of notation and conflate a map $\gamma:\bT\rightarrow\R^n$ with the image $\gamma(\bT)$ and, for example, write $\what\gamma$ instead of $\what{\gamma(\bT)}$ for the polynomial hull of the set ${\gamma(\bT)}$.

We recall the definition of the Hausdorff metric.  Let $X$ be a metric space with metric $d$, and let 
$\mathcal K$ denote the collection of all nonempty compact subsets of $X$. For $A,B\in\mathcal K$, the \emph{Hausdorff distance} $d_H(A,B)$ between them is defined to be the number
$$d_H(A,B)=\max
\{\max_{a\in A}\,\min_{b\in B}  \, d(a,b),\  \max_{b\in B}\, \min_{a\in A}\,  d(a,b)\}.$$
We will use the following standard compactness result.  An outline of the proof can be found in \cite[Section~45, Exercise~7]{Munkres:2000}.

\begin{thm}
\label{compactness}
If the metric space $X$ is compact, then the collection of all nonempty compact subsets of $X$ is a compact space with respect to the Hausdorff metric.
\end{thm}

%
%

\section{{\bf {Main Results}}.}\label{results}
\setcounter{newsection}{3}
\setcounter{thm}{0}

The following is our principal result regarding convergence of hulls:

\begin{thm}
\label{5.23.20.i}
If $\sqgamma$ is a sequence of rectifiable simple closed curves in $\cn$ that converges uniformly to a simple closed curve $\gamma$, and if each $\gamma_k$ has nontrivial polynomial hull, then $\gamma$ has nontrivial polynomial hull. If, in addition, the limit curve $\gamma$ is rectifiable, then the sequence $\sqhatgamma$ of polynomial hulls converges in the Hausdorff metric to $\what\gamma$.
 \end{thm}
 
Note that the first half of this theorem asserts something stronger than that the set of polynomially convex, rectifiable simple closed curves is open in the set of all rectifiable simple closed curves in the topology of uniform convergence: Given a polynomially convex simple closed curve $\gamma$, \emph{whether rectifiable or not}, there is an $\ve>0$ such that every rectifiable simple closed curve $\sigma$ satisfying $\|\gamma - \sigma \|_\bT<\ve$ is polynomially convex.
Since the topology of uniform convergence is weaker than the topology of $\bvt$,
the first half of Theorem~\ref{5.23.20.i} also implies an openness result in the space $\bvt$:
 
\begin{cor}
\label{5.25.20.v}
The set of polynomially convex, rectifiable simple closed curves is open in the set of all rectifiable simple closed curves in the topology of the Banach space $\bvt$.
\end{cor}

In \cite[Remark 1.6]{Nemirovski} Nemirovski presents an argument which, mutatis mutandis, shows that if $\gamma$ is a polynomially convex, \emph{rectifiable} simple closed curve, then every rectifiable simple closed curve sufficiently near $\gamma$ in the uniform norm is also polynomially convex (and hence, in particular, yields Corollary~\ref{5.25.20.v}). His argument, however, cannot be applied to the case of non-rectifiable simple closed curves $\gamma$.

In the context of Theorem~\ref{5.23.20.i}, it  is not claimed,  nor need it be true, that the sequence $\sqhatgamma$ of polynomial hulls converges to $\what\gamma$ when $\gamma $ is not rectifiable. See Example~\ref{example}  Note though that there is no requirement in Theorem~\ref{5.23.20.i} that the lengths of the rectifiable curves $\gamma_k$ be uniformly bounded.

It is clear that if a sequence $\sqgamma$ of simple closed curves in $\cn$ converges uniformly to a simple closed curve $\gamma$, then (regarded as a sequence of closed subsets of $\cn$) the sequence $\sqgamma$ converges to $\gamma$ in the Hausdorff metric.  The converse, however, does not hold.  We will show by simple examples that both halves of Theorem~\ref{5.23.20.i} become false if the hypothesis that $\sqgamma$ converges uniformly is replaced by the weaker hypothesis that $\sqgamma$ converges in the Hausdorff metric.

Theorem~\ref{5.23.20.i} is a result on simple closed curves.  We will present various examples illustrating difficulties that ensue when one tries to deal with more general kinds of sets.

Regarding density of polynomially convex curves, we will show that every rectifiable simple closed curve can be approximated in total variation norm by a polynomially convex, rectifiable simple closed curve that coincides with the original curve except on an arbitrarily small segment:

\begin{thm}
\label{5.25.20.vi}
Given a rectifiable simple closed curve $\gamma$ in $\cn$, $n\geq 2$, given $\ve>0$, and given an open ball $B$ of $\cn$ that intersects $\gamma$, there is a rectifiable simple closed curve $\gamma_a$ that is polynomially convex and satisfies $\|\gamma-\gamma_a\|_{\rm bv}<  \ve$ and 
$\gamma\setminus B=\gamma_a\setminus B$.
\end{thm}

\begin{cor}
\label{5.25.20.viiii} 
In $\cn,\ n\geq 2$, the set of polynomially convex, rectifiable simple closed curves is dense in the space $\bvt$ of all rectifiable closed curves with the total variation norm.
\end{cor}

This corollary is immediate from Theorem~\ref{5.25.20.vi} and the following result which we will prove regarding density of embeddings in the space of continuous maps of bounded 
variation.

\begin{thm}\label{maintheorem}
If $J$ is either a closed interval or a circle, and if $n\geq 3$, then the set of injective maps in $\bvjr$ is dense in $\bvjr$.
\end{thm}

Recall that an analogue of this result in the setting of smooth manifolds is a standard result in differential topology: {\it{ If $M$ is a compact manifold of dimension $d$ and of smoothness class $\sC^s,\ 1\leq s\leq\infty$, then the set of embeddings of class $\sC^s$ of $M$ into $\bR^k$ is dense in the space of all maps of class of $\sC^s$, in the $\sC^s$ topology, provided $k\geq2d+1$. }} This result is 
given in \cite[Chapter~2, Proposition~1.0]{Hirsch:1976}.  As $\sC^1(J)$ is not dense in $\bvjr$, Theorem~\ref{maintheorem} does not follow from this result.
For the same reason, the density of polynomially convex, rectifiable simple closed curves
in the space of all rectifiable simple closed curves with the total variation norm does not follow from the density of polynomially convex, $\sC^1$-smooth, simple closed curves
with the $\sC^1$ norm.  (That $\sC^1(J)$ is not dense in $\bvjr$ is easily seen by noting that $\sC^1(J)$ is contained in the set ${\rm AC}(J)$ of absolutely continuous functions and verifying that ${\rm AC}(J)$ is closed in $\bvjr$.  In fact, the closure of $\sC^1(J)$ in $\bvj$ is exactly ${\rm AC}(J)$ as can be seen by considering scalar-valued functions and noting that the map that sends each function of bounded variation on $J$ to the corresponding regular Borel measure is an isometry and sends $\sC^1(J)$ to the set $\sC(J)$ of continuous functions and ${\rm AC}(J)$ to $L^1(J)$.)

Note that the analogy that holds between rectifiable embeddings and smooth embeddings with regard to density does not carry over to openness; every rectifiable embedding can be modified by an arbitrarily small amount in total variation norm so as to become constant on a small interval.

A second corollary of Theorem~\ref{5.25.20.vi} concerns arcs:
\begin{cor}
\label{6.22.20.i}
A rectifiable arc in $\cn$, $n\geq 2$, is contained in a polynomially convex, rectifiable simple closed  curve, which can be chosen to lie in an arbitrarily small neighborhood of the given arc.
\end{cor}

Finally we will establish the following analogue of Theorem~\ref{5.25.20.vi} for smooth curves.
Here we denote by $d_{\sC^s}(\gamma, \gamma_a)$ the distance from $\gamma$ to $\gamma_a$ in $\sC^s(\bT)$.  This is of course given by a norm when $1\leq s<\infty$, but not when $s=\infty$.

\begin{thm}
\label{8.23.20.i}
Given a simple closed curve $\gamma$ of class $\sC^s$, $1\leq s\leq \infty$, in $\cn$, $n\geq 2$, given $\ve>0$, and given an open ball $B$ of $\cn$ that intersects $\gamma$, there is a simple closed curve $\gamma_a$ of class $\sC^s$ that is polynomially convex and that satisfies $d_{\sC^s}(\gamma, \gamma_a)<  \ve$ and 
$\gamma\setminus B=\gamma_a\setminus B$.
\end{thm}

The case of Theorem~\ref{8.23.20.i} in which $s=\infty$ follows immediately from \cite[Theorem~1.4]{Arosio-Wold:2019}.  However, we will give a single proof of Theorem~\ref{8.23.20.i} that applies in all cases and is much simpler than the proof of \cite[Theorem~1.4]{Arosio-Wold:2019}.   
The analogues of Corollaries~\ref{5.25.20.viiii} and~\ref{6.22.20.i} for smooth curves follow from Theorem~\ref{8.23.20.i}

\begin{cor}
\label{smooth-dense} 
In $\cn,\ n\geq 2$, the set of polynomially convex, simple closed curves of class $\sC^s$, $1\leq s\leq \infty$, is dense in the space $\sC^s(\T)$ of all closed curves of class $\sC^s$.
\end{cor}

\begin{cor}
\label{smooth-curve}
An arc of class $\sC^s$, $1\leq s\leq \infty$, in $\cn$, $n\geq 2$, is contained in a polynomially convex, simple closed curve of class $\sC^s$, which can be chosen to lie in an arbitrarily small neighborhood of the given arc.
\end{cor}

%
%

\section{{\bf {Convergence of Hulls}}.}\label{convergence}
\setcounter{newsection}{4}
\setcounter{thm}{0}

In this section we prove Theorem~\ref{5.23.20.i}  
and present related results and examples.
We will use the following three lemmas whose proofs we defer for the moment.

\begin{lem}
\label{6.20.20.i}
If  $\{X_k\}_{k=1,2,\dots}$ is a sequence of compact sets in $\cn$ such that the sequences $\{X_k\}_{k=1,2,\dots}$ and $\{\what X_k\}_{k=1,2,\dots}$ each converge in the Hausdorff metric, then $\lim \what X_k \subset \what{\lim  X_k}.$
\end{lem}

\begin{lem}\label{zeros}
Let $\sqgamma$ be a sequence of rectifiable simple closed curves each with nontrivial polynomial hull, and suppose the sequence $\sqgamma$ converges uniformly to a simple closed curve $\gamma$.  Let $f$ be a smooth $\bC$-valued function on $\cn$ whose restriction to $\gamma$ is
zero-free and has no continuous logarithm on $\gamma$.  Then $f$ has a zero on each subsequential limit of $\sqhatgamma$ in the Hausdorff metric.
\end{lem}

\begin{lem}
\label{5.25.20.i}
If $\gamma$ is a compact set that is contained in a connected set of finite length in $\cn$, and if $p$ is a point of $\what \gamma\setminus \gamma$, then there is a function $f$ holomorphic on a neighborhood of $\what \gamma$ that vanishes at $p$ and at no other point of $\wg$. 
\end{lem}

Note that because the function $f$ in Lemma~\ref{5.25.20.i} has a zero on the variety $\what \gamma\setminus \gamma$, the argument principle implies that $f$ has no continuous logarithm on $\gamma$.

\medskip

\noindent{\bf Proof of Theorem~\ref{5.23.20.i}}
The sequence of hulls $\sqhatgamma$ is contained in a compact subset of $\cn$, say a large closed ball.  Therefore, by Theorem~\ref{compactness}, this sequence has subsequential limits in the Hausdorff metric.  By Lemma~\ref{6.20.20.i}, each subsequential limit is contained in $\what\gamma$.  Thus to show that $\gamma$ has nontrivial polynomial hull, it suffices to show that some subsequential limit is not contained in $\gamma$.  Furthermore, to establish that $\sqhatgamma$ converges to $\what\gamma$ when $\gamma$ is rectifiable, it suffices to show that in that situation, every subsequential limit of $\sqhatgamma$ contains $\what\gamma$.

Consider a subsequence of the sequence $\sqhatgamma$ that converges in the Hausdorff metric, say to the set $Y$.  For notational convenience, we assume the subsequence to be the entire sequence itself.  Choose a continuous zero-free $\bC$-valued function on $\gamma$ that has no continuous logarithm.  Extend the function to a continuous function on $\cn$.  Finally by taking a sufficiently good approximation to that extension by a smooth function, obtain a smooth function $f$ whose restriction to $\gamma$ is zero-free and has no continuous logarithm on $\gamma$.
By Lemma~\ref{zeros}, $f$ has a zero on $Y$.  Hence $Y$, a subset of $\what\gamma$, is not contained in $\gamma$.  Thus the polynomial hull of $\gamma$ is nontrivial. 

When the limit curve $\gamma$ is rectifiable, Lemma~\ref{5.25.20.i} and the remark following it show that for any point $p$ of $\what \gamma\setminus \gamma$, the function $f$ can be chosen so that its only zero on $\what\gamma$ is at $p$.  Since $Y$ is contained in $\what\gamma$, and we have shown that $f$ must have a zero on $Y$, the point $p$ must belong to $Y$.  Consequently, in case $\gamma$ is rectifiable, the set $Y$ must contain $\what\gamma$.
\medskip

It remains to prove the lemmas.

\medskip

\noindent{\bf Proof of Lemma~\ref{6.20.20.i}}
This is perhaps well known, but we include the easy proof for completeness.
Set $X=\lim X_k$.
Let $P$ be an arbitrary polynomial on $\cn$, and let $\ve>0$ be arbitrary.  Then there exists $N$ such that for all $k\geq N$ we have
$$\|P\|_{X_k} \leq \|P\|_X +\ve.$$
This inequality continues to hold with $X_k$ replaced by $\what X_k$.  Since $\lim \what X_k$ is contained in the closure of the set $\bigcup_{k=N}^\infty \what X_k$, it follows that $\|P\|_{\lim \what X_k} \leq \|P\|_X +\ve$.  Consequently, $\|P\|_{\lim \what X_k} \leq \|P\|_X$.  Therefore, $\lim \what X_k \subset \what X.$

\medskip

\noindent{\bf Proof of Lemma~\ref{zeros}}
Consider a subsequence of the sequence $\sqhatgamma$ that converges in the Hausdorff metric.  For notational convenience, we assume the subsequence to be the entire sequence itself.  

Let  $U$ be the neighborhood of $\gamma$ given by $U=\{z\in \cn : f(z)\neq 0\}$. 
Because $\sqgamma$ converges uniformly to $\gamma$, for large $k$ the curve 
$\gamma_k$ is homotopic to $\gamma$ in $U$.
Accordingly, $f$ has no continuous logarithm on $\gamma_k$ for large $k$. 

The one-form $\omega$ defined locally to be $d\log f$ is well defined, smooth, and closed on $U$.  For $k$ large, $\int_{\gamma_k}\, \omega \neq 0$ since
$f$ has no continuous logarithm on $\gamma_k$. 

By hypothesis $\gamma_k$ is not polynomially convex, so the polynomial hull $\what\gamma_k$ is the union of $\gamma_k$ and a bounded purely one-dimensional variety $V_k$. (See \cite[Th. 3.1.1]{Stout:2007}.)
If $V_k$ were contained in $U$, then Stokes's Theorem would yield $\int_{\gamma_k}\, \omega=\int_{V_k}\, d\omega = 0$.
Since $\int_{\gamma_k}\, \omega \neq 0$ for $k$ large, we obtain that for $k$ large, the variety $V_k$ cannot be a subset of $U$.  Therefore, for $k$ large, the set $V_k$, and hence the set $\what\gamma_k$, meets the zero set of $f$.  Consequently, the limit set $\lim \what\gamma_k$ must also meet the zero set of $f$.
\medskip

\noindent{\bf Proof of Lemma~\ref{5.25.20.i}} If $\gamma$ is a set in the complex plane, the result is evident.

We will treat first the result in $\c2$ and then reduce the case of sets in $\cn$, $n>2$, to the case of sets in $\bC^2$.

Consider then a compact set $\gamma$ that is contained in a connected set of finite length in $\c2$. Let $V=\wg\setminus\gamma$.  Then $V$ is a purely one-dimensional analytic subvariety of $\cn\setminus\gamma$ (see, for instance, \cite[Theorem~3.1.1]{Stout:2007}).

Denote by $\sO$ the sheaf of germs of holomorphic functions on $\c2$ and by $\sO^*$ the sheaf of germs of zero-free holomorphic functions on the same $\c2$.
With the map $\sO\rightarrow \sO^*$ the map given by $f\mapsto e^{2\pi i f}$, there is the exact sequence of sheaves
$$
0\rightarrow\bZ\rightarrow\sO\rightarrow\sO^*\rightarrow 0.
$$
The associated cohomology sequence on $\what\gamma$
contains the segment
\begin{equation*}
\cdots\rightarrow H^1(\wg;\sO)\rightarrow H^1(\wg,\sO^*)\rightarrow H^2(\wg;\bZ)\rightarrow\cdots.
\end{equation*}
In this, $H^1(\wg;\sO)$ is the zero group because $\wg$ is the intersection of a decreasing sequence of Stein domains and because cohomology is continuous. 
We also have that $ H^2(\wg;\bZ)=0$, because $\wg$ is a polynomially convex set in $\c2$. (It is a well-known result of Andrew Browder  \cite{Browder:1961} that for a compact polynomially convex set $K$ in $\cn$, the cohomology groups $H^j(K;\bC)$ vanish for $j\geq n$.  It was observed in \cite{Duchamp-Stout:1981} that with integral coefficients, the analogous vanishing theorem is also correct. Vanishing theorems of this kind have been discussed in some detail in \cite{Stout:2007}.  The specific theorem we invoke here is \cite[Corollary~2.3.6]{Stout:2007}.)

Granted that $H^1(\wg;\sO)$ and $H^2(\wg;\bZ)$ both vanish
it follows that $H^1(\wg;\sO^*)=0$. 

Since the set $\gamma$ is of zero two-dimensional Hausdorff measure, it is rationally convex.
Thus there is a polynomial $\vp$ on $\c2$  that vanishes at $p$ but is zero-free on 
$\gamma$. 
This polynomial will be identically zero on no branch of  $V$, for otherwise its zero locus would meet $\gamma$.$^($\footnote{It should perhaps be observed that granted only that $\gamma$ is contained in a connected set of finite length, the variety $V$ is not assured to be irreducible; indeed, it could have infinitely many topological components.  However, each \emph{germ} of an analytic variety has only finitely many irreducible branches.}$^)$

As the intersection of the  zero locus of $\vp$ with each irreducible branch of $V$ is a (possibly empty) discrete set, there is an open ball $B$ in $\c2$ centered at $p$ such that $\overline B$ is disjoint from $\gamma$ and such that on 
$\overline B\cap V$ the polynomial $\vp$ vanishes only at $p$. 
 Let $D$ be a neighborhood in $\bC^2$ of the compact set $\what\gamma\setminus B$  that is disjoint from the set $\overline B\cap \vp^{-1}(0)$. We now have a set of Cousin~II data on the open set $D\cup B$: Take $\vp$ on $B$ and the function identically one on $D$. Because $H^1(\wg;\sO^*)=0$, this set of Cousin data is solvable on some neighborhood $\Omega$ of $\wg$,$^($\footnote{
A word of explanation may be in order here.
Let $\{\Delta_j\}_{j=1,2,\dots}$ be a decreasing sequence of compact neighborhoods of the set $\what\gamma$ with $\cap_{j} \Delta_j=\wg$ and with $\Delta_1\subset D\cup B$. Thus, $\wg$ is the {\it{inverse}} limit of the
$\Delta_j$ with the inclusion maps $\iota_{j,k}:\Delta_j\rightarrow \Delta_k$ for $j>k$, and the cohomology group $H^1(\wg;\sO^*)$ is the {\it{direct}} limit of the system $H^1(\Delta_j;\sO^*)$ with the induced maps $\iota^*_{j,k}:H^1(\Delta_k;\sO^*)\rightarrow H^1(\Delta_j;\sO^*)$ for $k>j$. This direct limit is $0$.

The set of Cousin II data we have constructed above on $D\cup B$ gives rise by restriction to a set of 
Cousin II data on each $\Delta_j$ and so for each $j$ a cohomology class $c_j\in H^1(\Delta_j;\sO^*)$. We have $\iota^*_{j,k}(c_k)=c_j$ for $k>j$. For each $j$,  there is the map $H^1(\Delta_j;\sO^*)\rightarrow H^1(\wg;\sO^*)$ which is the zero map. Consequently, for sufficiently large $j$, the cohomology class $c_j$ is zero. This means that our Cousin II problem is solvable on $\Delta_k$ for large $k$.}$^)$  so there is a function defined and holomorphic on a neighborhood of $\wg$ whose zero set meets $\wg$ only at $p$.  

To deduce the $\cn$ version of the result from the $\c2$ version, we proceed by projection.

Thus, let $\gamma$ be a compact set that is contained in a connected set of finite length in $\cn$, and let $p$ be a point in $V=\wg\setminus\gamma$. 
Exactly as in the case when $\gamma$ was in $\c2$, we can obtain a polynomial 
$\vp_1$ on $\cn$  that vanishes at $p$ but is zero-free on 
$\gamma$, and this polynomial will be identically zero on no branch of  $V$.

The set $\vp_1^{-1}(0)\cap \wg$ is polynomially convex and has $p$ as an isolated point. Consequently, there is a polynomial $\vp_2$ that vanishes at $p$ and at no other point of $\vp_1^{-1}(0)\cap \wg$. Let $\Phi=(\vp_1,\vp_2):\cn\rightarrow \c2$. The map $\Phi$ carries $\gamma$ to a compact set $\sigma$ in $\c2$  that is contained in a connected set of finite length, and it carries $\wg$ into, though perhaps not onto, the hull $\what\sigma$. Moreover, $\Phi(p)\notin\sigma$.

The version of our lemma already proved in $\c2$ provides a function $f$ holomorphic on a neighborhood of $\what\sigma$ whose zero set meets $\what\sigma$ only at the origin. The composition $f\circ\Phi$ is holomorphic on a neighborhood of $\what\gamma$ and vanishes at the point $p$ and at no other point of $\wg$.

The lemma is proved.
\medskip

With Lemmas~\ref{6.20.20.i}--\ref{5.25.20.i} established, the proof of Theorem~\ref{5.23.20.i} is complete.

As mentioned in Section~\ref{results},
in the context of Theorem~\ref{5.23.20.i}, the sequence $\sqhatgamma$ of polynomial hulls can fail to converge to $\what\gamma$ when $\gamma $ is not rectifiable. An example of this phenomenon is the following.

\smallskip

\begin{example}\label{example}
{\rm Let $C$ be a simple closed curve of positive area in the plane$^($\footnote{Such curves were constructed by Osgood \cite{Osgood:1903}.}$^)$,
and let $F=(f_1,f_2,f_3)$ be a continuous injective map from the Riemann sphere into $\bC^3$ that is holomorphic on the interior $D_i$ of $C$ and also on the exterior $D_e$ of $C$. Such maps were considered by Wermer \cite{Wermer:1955} as follows.
Set  
$$f_1(z)=\int_C\frac{d\z\wedge d\overline \z}{\z-z}\ \ 
{\rm{and}}\ \ 
f_2(z)=zf_1(z).$$
Then choose a point $z_0$ in $\bC\setminus C$ at
which $f_1(z)$ does not vanish and define $f_3$ by
$$f_3(z)=\frac{f_1(z)-f_1(z_0)}{z-z_0}.$$
With these three functions, the map $F=(f_1,f_2,f_3)$ carries the Riemann sphere continuously and injectively into $\bC^3$ and is holomorphic off the curve $C$.
 
Let $\vp$ be a conformal map from the unit disc $\bU$ in $\bC$ onto $D_i$. The map $\vp$ extends to a homeomorphism of $\overline \bU$ onto 
$\overline D_i$. The existence of this extension is a theorem of Carath\'eodory which is given in \cite[p. 13]{Garnett-Marshall:2005}. 

For $k=2,3,\dots$ let $\beta_k$ be the boundary of the disc $\Delta_k=\{z\in\bC:|z|\leq 1-\frac{1}{k}\}$. 
Let $\gamma_k$ be  $F( \vp(\beta_k))$. Each $\gamma_k$ is a real-analytic simple closed curve, and $\gamma_k$ tends uniformly to the curve $\gamma$ obtained from $F\circ \vp$ by restricting to the boundary of $\bU$.

We do not know what $\what\gamma$ is, but it does contain the sets $F(D_i)$ and $F(D_e)$.
For each $k$, it is clear that $\what\gamma_k\supset F(\vp(\Delta_k))$.  The set $F(\vp(\Delta_k))\setminus \gamma_k$ is an analytic subvariety of $\cn\setminus\gamma_k$. (This can be seen by showing that the derivative of the map $F$ is nowhere vanishing on $\bC\setminus C$, and hence $F(\bC\setminus C)$ is, in fact, a complex manifold.  Alternatively, it is a very minor case of  Remmert's proper mapping theorem \cite[Theorem~N1]{Gunning:1990}). 
 Since the variety 
$\what\gamma_k\setminus\gamma_k$ is irreducible \cite[Theorem~4.5.5]{Stout:2007}, $\what\gamma_k\setminus\gamma_k$ can be no larger than $F(\vp(\Delta_k))\setminus \gamma_k$.  Therefore, $\what\gamma_k=F(\vp(\Delta_k))$.

The sequence of sets $\sqhatgamma=\{F(\vp(\Delta_k))\}_{k=1,2,\ldots}$ is increasing and has union $F(\vp(\bU))$.  It follows that $\what\gamma_k\rightarrow F(\overline D_i)$.
Consequently, $\sqhatgamma$ does not converge to $\what\gamma$.}

\end{example}

As mentioned in Section~\ref{results}, both halves of Theorem~\ref{5.23.20.i} become false if the hypothesis that $\sqgamma$ converges to $\gamma$  uniformly is replaced by the weaker hypothesis  that $\sqgamma$ converges to $\gamma$ in the Hausdorff metric.  This is demonstrated by the next two examples.

\begin{example}{\rm
Let $\gamma$ be the unit circle in the plane.  For each $k$, let $\lambda_k$ be the 
arc on the circle $\{ z: |z|= 1 + 1/2k\}$ consisting of those points whose argument lies in the interval $[\pi/2k, 2\pi - \pi/2k]$. Let $\sqgamma$ be a sequence of smooth simple closed curves in the plane such that, for each $k$, the arc $\lambda_k$ is contained in the bounded component of the complement of $\gamma_k$, and such that $\gamma_k$ is contained in the slit annulus obtained from the annulus $\{z: 1<|z|< 1+1/k\}$ by deleting the positive real axis. (See Figure~1.)  Then $\gamma_k\rightarrow \gamma$ in the Hausdorff metric, but $\what\gamma_k\rightarrow \gamma\neq \what\gamma=\overline\bU$.}
\end{example}

\vskip 74 pt
\hskip2.73 truein
\psset{linewidth=0.1pt}
\psset{xunit=0.16cm}
\psset{yunit=0.16cm}
\psset{runit=0.16cm}
\pscircle[linecolor=black,linewidth=1.1pt,](0,0){13.75}
\psarc[linecolor=black,linewidth=1.1pt,](0,0){15.5}{5}{355}
\psarc[linecolor=black,linewidth=1.1pt,](0,0){16.5}{5}{355}
\psarc[linecolor=black,linewidth=1.1pt,](15.95,1.4){0.5}{180}{360}
\psarc[linecolor=black,linewidth=1.1pt,](15.95,-1.4){0.5}{0}{180}

\vskip -55pt
\hskip3.2 truein \hskip -16pt
$\gamma$
\vskip -7pt
\hskip3.25 truein \hskip 41pt $\gamma_k$
     
\vskip 120 truept
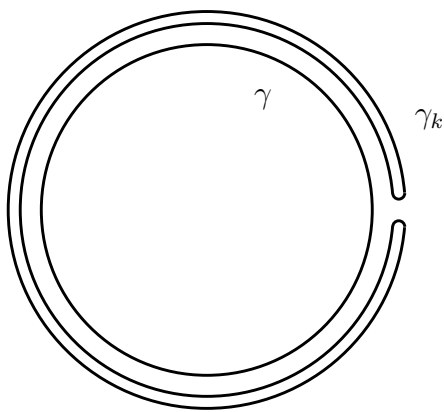
\centerline{Figure 1}
\vskip 40 pt

\psset{linewidth=0.8pt}   

\begin{example}{\rm
Let $\sigma$ and $\sigma_k$ be, respectively, the images of the simple closed curves $\gamma$ and $\gamma_k$ of the previous example under the map of $\bC$ into $\bC^2$ given by $z\mapsto (z, 1/z)$. Then $\sqsigma$ converges to $\sigma$ in the Hausdorff metric, but $\sigma$ is polynomially convex while each $\sigma_k$ boundes an analytic disc on the analytic variety $\{(z, 1/z): z\in \bC\setminus \{0\}\}$
and hence has nontrivial polynomial hull.}
\end{example}

Theorem~\ref{5.23.20.i} is a result about simple closed curves, or in other words, about topological embeddings of a circle into $\cn$.  The following two examples show that both halves of Theorem~3.1 fail when the circle is replaced by a more general compact, connected space, even under the additional hypothesis that there is a uniform bound on the lengths of the embedded sets.  

\begin{example}{\rm
Let $K$ be the subspace of the plane that is the union of the unit circle $b\bU$ and a countable collection of circles $C_k$, $k=1,2,\ldots$, with each $C_k$ externally tangent to $b\bU$ at the point 
$e^{i\pi/k}$ and of radius $|e^{i\pi/k} - e^{i\pi/(k+1)}|/4$.  Note that $K$ is a compact, connected space.  Let $p_k=(e^{i\pi/k},e^{-i\pi/k})$.  
For each $k=1,2,\ldots$, choose circles $G_k$ and $E_k$ in $\bC^2$ of radius $|p_k - p_{k+1}|/4$ that intersect the circle $ \{(e^{i\vt}, e^{-i\vt}):0\leq \vt\leq  2\pi\}$ only in the point $p_k$ and such that $G_k$ is contained in the plane 
 $\{ (z, \overline z) : z \in \bC\}$ and $E_k$ is contained in the plane $\{ p_k + (z,0) : z\in \bC\}$.  Let 
$$ X=\{(e^{i\vt}, e^{-i\vt}):0\leq \vt\leq  2\pi\} \cup \bigcup_{k=1}^\infty G_k.$$
Set $X_k=(X\setminus G_k) \cup E_k$ so that $X_k$ is the set obtained from $X$ by removing the circle $G_k$ and replacing it with the circle $E_k$.  
Let $\rho:K\rightarrow X$ and $\rho_k:K\rightarrow X_k$ be the obvious homeomorphisms.
Then the $X_k$ all have the same finite length, each $X_k$ has nontrivial polynomial hull, $\rho_k\rightarrow \rho$ uniformly, 
but $X$ lies in the plane $\{ (z, \overline z) : z \in \bC\}$ and hence is polynomially convex.
}
\end{example}

\begin{example}{\rm
Let $K$ be the union of two circles meeting in a single point given by
$$K=\{\, z: |z|=1\} \cup \{\, 2+z: |z|=1\}.$$
Note that $K$ is compact and connected.
Set 
$$E_k=\{ (2+z, 1/k): |z|=1\},$$
set
$$X_k=\{ (z, \overline z/k): |z|=1\} \cup E_k,$$
and set
$$X=\{ (z,0): |z|=1\} \cup \{ (2+z, 0): |z|=1\}.$$
Let $\rho:K\rightarrow X$ and $\rho_k:K\rightarrow X_k$ be the obvious homeomorphisms.
Then the lengths of the $X_k$ are bounded by $2(\sqrt 2 +1)\, \pi$ and $\rho_k\rightarrow \rho$ uniformly.
By applying Kallin's lemma \cite{Kallin:1965} (or see \cite[Theorem~1.6.19]{Stout:2007}), one obtains that 
$$\what X_k=\{ (z, \overline z/k): |z|=1\} \cup \{ (2+z, 1/k): |z|\leq1\}.$$
Thus
\begin{equation*}\begin{split}\lim \what X_k &= \{ (z,0): |z|=1\} \cup \{ (2+z, 0): |z|\leq 1\}\\ &\subsetneq \{ (z,0): |z|\leq1\} \cup \{ (2+z, 0): |z|\leq 1\} = \what X.\\
\end{split}
\end{equation*}
}
\end{example}
\medskip

%
%

\section{Density of Rectifiable Embeddings.}\label{embeddings}
\setcounter{newsection}{5}
\setcounter{thm}{0}

This section is devoted to proving Theorem~\ref{maintheorem} on the density of rectifiable embeddings in the space of continuous maps of bounded variation.
The proof depends on two lemmas, the first a simple result in geometric measure theory.  

We denote the $k$-dimensional Hausdorff measure of a set $E$ in $\bR^n$ by 
${\scr H}^k(E)$.  By the length of a rectifiable curve $\gamma$, we mean the total variation of $\gamma$.  Note that in case $\gamma$ is not injective, the length of $\gamma$ may well exceed ${\scr H}^1(\gamma(J))$.

\begin{lem} If $J$ is either a closed interval or a circle, and if $\gamma:J\rightarrow \R^n$  is a rectifiable curve of length $l$, then ${\scr H}^2(\gamma\times\gamma)\leq (\pi/2)l^2$. In particular, $\gamma\times\gamma$ has finite 2-dimensional Hausdorff measure.
\end{lem}

\smallskip

\noindent{\bf{Proof.}} Fix $\vep >0$.  The lemma will be established once we show that there exists a countable collection of sets $A_1, A_2, \ldots$ that covers $\gamma\times \gamma$ with each set $A_j$ of diameter $\delta(A_j)< \vep$ and such that 
\begin{equation*}
(\pi/4) \sum_j \delta(A_j)^2 \leq (\pi/2) l^2.
\end{equation*}

Choose $m\in \Z_+$ large enough that $l/m< \vep/\sqrt 2$.  Partition $J$ into $m$ subintervals $J_1,\dots,J_m$ such that each of the restrictions $\gamma_j=\gamma|J_j$ has length
 $l/m$.  Then $\gamma\times \gamma = \bigcup_{j,k=1}^m \gamma_j\times \gamma _k$.  A trivial computation shows that $\delta(\gamma_j \times \gamma_k) \leq \sqrt 2 (l/m)$ for each $j$ and $k$ so that
\begin{eqnarray*}
\sum_{j,k=1}^m \delta(\gamma_j\times\gamma_k)^2 & \leq& \sum_{j,k=1}^m 2(l/m)^2\cr
&=& 2 m^2 (l/m)^2\cr
&=& 2 l^2.\cr
\end{eqnarray*}
Thus
\begin{equation*}
(\pi/4) \sum_{j,k=1}^m \delta(\gamma_j\times \gamma_k)^2 \leq (\pi/2) l^2,
\end{equation*}
and the lemma is proved.

\begin{lem}Let $J$ be either a closed interval or a circle.    Let $\gamma : J\rightarrow \R^n$, $n\geq4$, be an injective continuous map of bounded variation, and let $\vep>0$.  Let $P: \R^n \rightarrow \R^{n-1}$ denote the projection onto the last $n-1$ coordinates.  Then there exists a linear operator $T:\R^n\rightarrow \R^{n-1}$ with $\|T-P\|<\vep$ such that $T\circ \gamma$ is injective.
\end{lem}

\noindent{\bf Proof.}
Let $\widetilde P$ denote the orthogonal projection of $\R^n$ onto $\{0\} \times \R^{n-1}$.
Given $v$ in the unit sphere $\bS^{n-1}$ in $\R^n$ with $v$ not in $\{0\} \times \R^{n-1}$, let $T_v:\R^n\rightarrow \R^n$ be the linear projection with range $\{0\} \times \R^{n-1}$ and null space the linear span of $v$.  It suffices to show that the set of vectors $v$ such that $T_v\circ \gamma$ is injective is dense in 
$\bS^{n-1}$, for if $v$ is sufficiently close to the vector $(1,0,\ldots, 0)$, then $\|T_v - \widetilde P\|<\vep$, and hence setting $T=P\circ T_v$ yields the lemma.

Let $\Delta$ denote the diagonal $\Delta= \{(x,y)\in \R^n\times \R^n: x=y\}$ in $\R^n\times\R^n$, and define $g:(\R^n\times \R^n)\setminus \Delta\rightarrow \bS^{n-1}$ by
\begin{equation*}
g(x,y)= \frac{x-y}{|x-y|}.
\end{equation*}
Then, for $v$ in $\bS^{n-1}\setminus\bigl(\{0\} \times \R^{n-1}\bigr)$, the map $T_v\circ\gamma$ is injective if and only if $v$ is not in $g\bigl((\gammasq)\setminus\Delta\bigr)$.  By the above lemma, $\gammasq$ has finite 2-dimensional Hausdorff measure.  Since $g$ is a smooth map on $(\R^n\times \R^n)\setminus \Delta$, it follows that $g$ maps compact subsets of $(\gammasq)\setminus\Delta$ to sets of finite 2-dimensional Hausdorff measure in $\bS^{n-1}$, and so, in particular, $g\bigl((\gammasq)\setminus\Delta\bigr)$ has 3-dimensional Hausdorff measure zero.  Thus $g\bigl((\gammasq)\setminus\Delta\bigr)$ has empty interior in $\bS^{n-1}$, and the lemma is proved.
\medskip

\noindent{\bf Proof of Theorem~\ref{maintheorem}}
Let $\gamma:J\rightarrow \R^n$ be a continuous map of bounded variation, and let $\vep>0$.
In case $J$ is the unit interval, let $\sigma:J\rightarrow \R^{1+n}$ be the graphing map $\sigma(x)=\bigl( x, \gamma(x)\bigr)$.  By the preceding lemma, there is a linear operator $T:\R^{1+n}\rightarrow \R^{n}$ such that $T\circ\sigma$ is injective and $\|T-P\|<\vep/\|\sigma\|_{\rm bv}$, where $P$ is the projection $\R^{1+n}\rightarrow\R^{n}$ onto the last $n$ coordinates.  Then
\begin{equation*}
\| (T\circ \sigma) - \gamma \|_J = \| (T\circ \sigma) - (P\circ \sigma)\|_J \leq \|T-P\|\, \| \sigma \|_J
\end{equation*}
and
\begin{equation*}
{\rm var} \bigl((T\circ \sigma) - \gamma \bigr)\leq {\rm var}\bigl((T\circ \sigma) - (P\circ \sigma)\bigr) \leq \|T-P\|\, {\rm var}\, \sigma 
\end{equation*}
so $\|(T\circ \sigma) -\gamma\|_{\rm bv}< \vep$.

The proof when $J$ is a circle is the same except that one needs to apply the lemma twice since in that case the graph of $\gamma$ lies in $\R^{2+n}$.

%
%

\section{Density of Polynomially Convex Simple Closed Curves.}\label{density}
\setcounter{newsection}{6}
\setcounter{thm}{0}

In this section we prove Theorems~\ref{5.25.20.vi} and~\ref{8.23.20.i} and Corollaries~\ref{6.22.20.i} and~\ref{smooth-curve}  The proofs are based on 
the following characterization of those rectifiable simple closed curves that are polynomially convex.

\begin{thm}
\label{5.25.20.ix}
The rectifiable simple closed curve $\gamma$ in $\cn$ is polynomially convex if and only if there is a holomorphic one-form $\alpha$ on $\cn$ such that $\int_\gamma\, \alpha \neq 0$.
\end{thm}
 
This result is given in \cite[p.194]{Stout:2007}. For the convenience of the reader, we recall its brief proof.  (Here
we denote by
$\sC(\gamma)$ the algebra of all continuous $\bC$--valued functions on $\gamma$, 
by $\sR(\gamma)$ the subalgebra of $\sC(\gamma)$ comprising those functions that can be approximated uniformly on $\gamma$ by rational functions holomorphic on a neighborhood of $\gamma$, and by 
$\sP(\gamma)$ the subalgebra of $\sC(\gamma)$ comprising those functions that can be approximated uniformly on $\gamma$ by holomorphic polynomials.)
 \medskip

\noindent{\bf Proof.} If $\gamma$ is not polynomially convex, its polynomial hull is $\gamma\cup V$  with $V$ an irreducible one-dimenional variety. Stokes' theorem, which is valid in this context \cite{Lawrence:1995},\cite[p.193]{Stout:2007}, gives 
$$\int_\gamma\, \alpha=\int_V\,d\alpha =0$$
for every holomorophic one-form $\alpha$ on $\cn$ because the holomorphic two-form $d\alpha$ vanishes on the one-dimensional variety $V$. 
On the other hand, if $\gamma$ is polynomially convex, then $\sP(\gamma)=\sC(\gamma)$.  (To see this, note that the rectifiability of $\gamma$ implies that $\sR(\gamma)=\sC(\gamma)$ by \cite[Theorem~1.6.7]{Stout:2007} for instance, and the polynomial convexity of $\gamma$ implies that $\sP(\gamma)=\sR(\gamma)$.)  It follows that there exists a holomorphic one-form $\alpha$ on $\cn$ such that $\int_\gamma\, \alpha\neq 0$. The theorem is proved.
\medskip

\noindent{\bf Proof of Theorem~\ref{5.25.20.vi}} 
Choose a point $p$ of $B \cap \gamma$, and let $B_p$ be an open ball centered at $p$ with radius less than $\ve/8$ and small enough that $B_p$ is contained in $B$, that the set $\gamma \cap B_p$ is contained in an arc $\lambda$ in $\gamma$ of length less than 
$\ve/8$.  Let $\Lambda$ be the component of $\gamma\setminus B_p$ that contains $\gamma\setminus \lambda$.  Note that $\Lambda$ is an arc and that its end points, which we will denote by 
$a$ and $b$, are in $bB_p$.

Given points $x$ and $y$, let $[x,y]$ denote the straight line segment from $x$ to $y$.  
Choose a point $c\in B_p$ that is not on the complex line through $a$ and $b$.  Let $\sigma$ denote the simple closed curve that is the boundary of the triangle with vertices $a$, $b$, and $c$ oriented so that $\sigma=[a,b]\cup[b,c]\cup[c,a]$.  

Introduce two rectifiable simple closed curves $\gamma+$ and $\gamma^-$ as sets by 
$$\gamma^+=\Lambda\cup [a,b]\ {\rm{and}}\ 
\gamma^-=\Lambda\cup [a,c] \cup[c,b].$$
As maps from the circle, define $\gamma^+$ and $\gamma^-$ to each coincide with the map $\gamma$ on the set $\gamma^{-1}(\Lambda)$ and to map $\gamma^{-1}(\gamma\setminus \Lambda)$ one-to-one onto $[a,b]$ and $[a,c]\cup[c,b]$, respectively, traversed in the direction that yields well-defined continuous maps.  

The simple closed curve $\sigma$ is polynomially convex, since it is contained in a totally real plane.  
Thus by Theorem~\ref{5.25.20.ix}, there is a holomorphic one-form $\alpha$ on $\cn$ such that $\int_{\sigma} \alpha \neq 0$.  For this $\alpha$ we have
\begin{equation*}
\label{5.25.20.nonzero}
\int_{\gamma^+}\!\alpha\, -\int_{\gamma^-}\!\alpha\,=\int_{\sigma}\alpha
 \neq 0,
\end{equation*}
so at least one of $\int_{\gamma^+}\alpha$ and $\int_{\gamma^-}\alpha$ is nonzero.
Consequently, by Theorem~\ref{5.25.20.ix}, at least one of $\gamma^+$ and $\gamma^-$ is polynomially convex.  

Observe that $\| \gamma - \gamma^\pm\|_{\bT}$ is bounded above by the sum of the length of $\lambda$ and the radius of $B_p$, and ${\rm{var}}(\gamma -\gamma^\pm)$ is bounded above by the sum of the length of $\lambda$ and twice the diameter of $B_p$.  Consequently, 
$\|\gamma -\gamma^\pm\|_{\rm bv}\leq 7\ve/8  < \ve$. 
Also $\gamma\setminus B=\gamma^\pm\setminus B$.

The theorem is proved.
\medskip

\medskip

\noindent{\bf Proof of Theorem~\ref{8.23.20.i}}
The map $\gamma:\bT\rightarrow \cn$ is of class $\sC^s$, is injective, and has nonvanshing derivative at each point of $\bT$. 

Given the ball $B$, fix a point $p\in B\cap\gamma$.  Without loss of generality $p=\gamma(1)$.  Let $\Delta$ denote the diagonal $\Delta= \{(z,w)\in \cn\times \cn: z=w\}$ in $\cn\times\cn$, and define $g:(\cn\times \cn)\setminus \Delta\rightarrow \bS^{2n-1}=\{ z\in \cn:|z|=1\}$ by
\begin{equation*}
g(z,w)= \frac{z-w}{|z-w|}.
\end{equation*}
Since $\gamma\times\gamma$ is a smooth 2-dimensional manifold and $g$ is a smooth map, the set $g\bigl((\gammasq)\setminus\Delta\bigr)$ has measure zero in the sphere 
$\bS^{2n-1}$.  Therefore, we can choose a unit vector $\bv$ not in $g\bigl((\gammasq)\setminus\Delta\bigr)$ and such that the real-linear span of $\bv$ and the tangent vector to $\gamma$ at $p$ is a totally real two-plane.
Define $G:\bT\times\bR\rightarrow\cn$ by
$$G(e^{i\vt},t)=\gamma(e^{i\vt})+t\bv.$$ 
By our choice of $\bv$, the map $G$ is injective and carries some neighborhood of $(1,0)$ in $\bT\times\bR$ onto a totally real manifold through $p$ in $\cn$.
By the local polynomial convexity of totally real manifolds in $\cn$, we can choose an interval $I$ in $\bT$ centered at the point $1$ and an $\eta>0$ such that every compact subset of $G(I\times[-\eta,\eta])$ is polynomially convex.
By choosing $I$ and $\eta$ small enough, we can also arrange to have $G(I\times[-\eta,\eta])\subset B$.

Choose a nonnegative function $\chi$ of class $\sC^\infty$ defined on $\bT$ such that the support of $\chi$ is a nonempty interval $I_0$ contained in the interior of the interval $I$ and such that $d_{\sC^s}(\chi,0)<\min\{\eta,\ve\}$.

Define maps $\gamma^+$ and $\gamma^-$ from $\bT$ to $\cn$
by
$$\gamma^+(e^{i\vt})=\gamma(e^{i\vt})+\chi(e^{i\vt})\bv\ \ {\rm{and}}\ \ \gamma^-(e^{i\vt})=\gamma(e^{i\vt})-\chi(e^{i\vt})\bv.$$
By our choice of $\bv$, these maps are both simple closed curves.  

The present argument now finishes along the lines of the previous proof:
Let $\sigma$ be the simple closed curve $\gamma^+(I_0)\cup\gamma^-(I_0)$, which is not smooth but is rectifiable.  Then $\sigma$ is polynomially convex since it is contained in the set $G(I\times [-\eta, \eta])$.  Thus there is a holomorphic one-form $\alpha$ on $\cn$ such that 
$$\int_{\gamma^+}\!\alpha\,-\int_{\gamma^-}\!\alpha=\int_\sigma\alpha\neq 0,$$
whence at least one of $\gamma^+$ and $\gamma^-$ is polynomially convex.

The curves $\gamma^\pm$ satisfy $d_{\sC^s}(\gamma, \gamma^\pm)=d_{\sC^s}(\chi,0)<\ve$
and $\gamma\setminus B=\gamma^\pm\setminus B$.

The theorem is proved.

\medskip

Corollary~\ref{6.22.20.i} is a consequence of Theorem~\ref{5.25.20.vi} and the following lemma.

\begin{lem}\label{8.21.20.i} 
If $\lambda$ is a rectifiable arc in $\bR^n,\ n\geq 2$, and $\Omega$ is a neighborhood of $\lambda$, then there is a rectifiable simple closed curve $\gamma$ that contains $\lambda$ and that is contained in $\Omega$.
\end{lem}

\noindent{\bf Proof.}
We will treat first the case $n\geq 3$, and then give a different argument for the case $n=2$.  Note that we need the result only for $n\geq 4$, so the argument for the case $n=2$ can be omitted.

Suppose then that $n\geq 3$.
We may assume without loss of generality that $\Omega$ is connected.
Let the end points of $\lambda$ be $p_0$ and $q_0$. Because the projective space of real lines in $\R^n$ through $p_0$ has dimension $n-1$ and $\lambda$ has finite length, there is a real line that passes through $p_0$ and is otherwise disjoint from $\lambda$. Let $p_1$ be a point on this line such that the straight line segment 
$[p_0,p_1]$ is contained in $\Omega$. Similarly, there is a real line through $q_0$ that is disjoint from $(\lambda\setminus\{q_0\}) \cup [p_0,p_1]$. Let $q_1$ be a point on this line such that the straight line segment 
$[q_0,q_1]$ is contained in $\Omega$.
Choose open Euclidean balls $B_p$ and $B_q$ centered at $p_1$ and $q_1$, respectively, such that the closures of $B_p$ and $B_q$ are disjoint and lie in $\Omega$, such that $B_p$ is disjoint from $\lambda \cup [q_0,q_1]$, and $B_q$ is disjoint from $\lambda\cup[p_0,p_1]$.  Choose points $p'$ and $q'$ in $B_p\setminus [p_0,p_1]$ and $B_q\setminus [q_0,q_1]$, respectively.  The set $\Omega\setminus \bigl(\lambda \cup [p_0,p_1] \cup [q_0,q_1]\bigr)$ is connected (because a connected manifold of real dimension greater than or equal to three cannot be disconnected by a subspace of topological dimension one \cite[Corollary~1,~p.~48]{Hurewicz-Wallman:1948}), so there is a rectifiable arc from $p'$ to $q'$ in 
$\Omega\setminus \bigl(\lambda \cup [p_0,p_1] \cup [q_0,q_1]\bigr)$.  By discarding initial and final segments of this arc, we can obtain an arc $\ell$ in 
$\Omega\setminus(\lambda \cup [p_0,p_1] \cup [q_0,q_1] \cup B_p \cup B_q)$ whose end points $p_2$ and $q_2$ lie on the boundary of $B_p$ and $B_q$, respectively.  Let $[p_1,p_2]$ denote the straight line segment from $p_1$ to $p_2$ and similarly with $p$ replaced by $q$.
Then $\lambda\cup [p_0,p_1] \cup [p_1,p_2] \cup \ell \cup [q_1,q_2] \cup [q_0, q_1]$  is a rectifiable simple closed curve in $\Omega$ that contains $\lambda$.
This completes the proof in the case $n\geq3$.

For the case $n=2$, identify $\R^2$ with $\bC$, and let $\bC^*$ denote the Riemann sphere.  We will use the following theorem which appears with proof as \cite[Theorem~2.1]{Church:1964} and seems to be due to Marie Torhorst \cite{Torhorst:1921} but forgotten.

\begin{thm}
If $G$ is a simply connected region on $\bC^*$ such that $bG$ is a nondegenerate Peano continuum, then each prime end is a single point.  Thus, if $f$ is any conformal homeomorphism of $\bU$ onto $G$, then $f$ can be extended to $\overline\bU$ to be continuous.
\end{thm}

Note that $\bC^*\setminus\lambda$ is a simply connected domain in $\bC^*$ whose boundary $\lambda$ is a nondegenerate Peano continuum.  Thus by the theorem just quoted, each conformal homeomorphism of $\bU$ onto $\bC^*\setminus \lambda$ extends to a continuous map of $\overline\bU$ onto $\bC^*$.  Choose such a map $f$ that takes the point $1$ to an end point $e_1$ of $\lambda$ and sends the point $-1$ to the other end point $e_{-1}$ of $\lambda$.  
Let $\alpha$ denote the point of $\bU$ such that $f(\alpha)=\infty$.  Then the function $g$ given by $g(z)=(z-\alpha)^2 f'(z)$ lies in the Hardy space $H^1(\bU)$ \cite[pp.~221-222]{Garnett-Marshall:2005}.

Choose $r$ such that $|\alpha| <r<1$ and such that $r$ is large enough that the set $f(\{z: r\leq |z|\leq 1\})$ is contained in $\Omega$. Let $D_1$ denote the disc with center $(1+r)/2$ and radius $(1-r)/2$, so that the boundary of $D_1$ is a circle contained in the annulus $\{z:r\leq |z|\leq 1\}$ passing through the points $r$ and $1$.  Let $D_{-1}$ denote the disc obtained by reflecting $D_1$ through the imaginary axis.  

Since the function $g$ is in $H^1(\bU)$, the function $g$ has a harmonic majorant on $\bU$, that is, there exists a harmonic function $u$ on $\bU$ such that $|g(z)|\leq u(z)$ for all $z\in\bU$ \cite[Theorem~2.12]{Duren:1970}.  Consequently, $f'$ has a harmonic majorant on $D_1$ and hence lies in $H^1(D_1)$.  It follows that the curve $f(bD_1)$ has finite length.  Similarly, the curve $f(bD_{-1})$ has finite length as well.

We conclude that if we let $\ell_{-1}$ denote a semicircular arc along the circle $b D_{-1}$ from $-1$ to $-r$, let $\ell$ denote a semicircular arc along the circle $\{z:|z|=r\}$ from $-r$ to $r$, and let $\ell_1$ denote a semicircular arc along the circle $bD_1$ from $r$ to $1$, then $f(\ell_{-1} \cup \ell \cup \ell_1) \cup \lambda$ is a rectifiable simple closed curve in $\Omega$ containing $\lambda$.

The lemma is proved.
\medskip

Corollary~\ref{smooth-curve} is a consequence of Theorem~\ref{8.23.20.i} and the following smooth analogue of the preceding lemma.

\begin{lem}\label{smooth-closing-up} 
If $\lambda$ is an arc of class $\sC^s$, $1\leq s\leq \infty$, in $\bR^n,\ n\geq 2$, and $\Omega$ is a neighborhood of $\lambda$, then there is a simple closed curve $\gamma$
of class $\sC^s$ that contains $\lambda$ and that is contained in $\Omega$.
\end{lem}

\noindent{\bf Proof.}
The outline of the proof is similar to that of the $n\geq3$ case of the previous lemma.  
We may assume without loss of generality that $\Omega$ is connected.
Throughout the proof, by \emph{smooth} we shall mean of class $\sC^s$.
Let the end points of $\lambda$ be $p_0$ and $q_0$.
In some smooth local coordinate system about the point $p_0$, the arc $\lambda$ is just a straight line segment ending at $p_0$.  Therefore, there is an arc $\lambda_p$ from $p_0$ to another point $p_1$ such that the union $\lambda \cup \lambda_p$ is a smooth arc contained in $\Omega$.
Similarly, there is an arc $\lambda_q$ from $q_0$ to another point $q_1$ such that the union $\lambda \cup \lambda_p \cup \lambda_q$ is also a smooth arc contained in 
$\Omega$.
Choose open Euclidean balls $B_p$ and $B_q$ centered at $p_1$ and $q_1$,  respectively, 
such that the closures of $B_p$ and $B_q$ are disjoint and lie in $\Omega$, such that $B_p$ is disjoint from $\lambda \cup\lambda_q$, and $B_q$ is disjoint from $\lambda\cup\lambda_p$.  Choose points $p'$ and $q'$ in $B_p\setminus \lambda_p$ and $B_q\setminus\lambda_q$, respectively.  The set 
$\Omega\setminus(\lambda \cup \lambda_p \cup \lambda_q)$ is connected. (When $n\geq 3$, this is immediate from dimensional considerations \cite[Corollary~1,~p.~48]{Hurewicz-Wallman:1948}.  To see that it holds also when $n=2$, first show that every connected neighborhood $U$ of an arc $J$ in the plane contains a connected neighborhood $V$ whose complement in the plane is connected.  Since $V$ is then homeomorphic to the plane, it is a standard fact that $V\setminus J$ is connected.  Connectedness of $U\setminus J$ follows.)
Therefore, there is a smooth arc $\ell$ from $p'$ to $q'$ in 
$\Omega\setminus(\lambda \cup \lambda_p \cup \lambda_q)$.  
By perturbing the radii of the open balls $B_p$ and $B_q$, we may assume that their boundaries intersect $\ell$ and $\lambda\cup\lambda_p\cup\lambda_q$ transversally.    Let $\lambda_*$ be the subarc of $\lambda\cup\lambda_p\cup\lambda_q$ that is disjoint from $B_p\cup B_q$ and connects a point of $bB_p$ to a point of $bB_q$.  
Define $\ell_*$ similarly but with $\lambda\cup\lambda_p\cup\lambda_q$ replaced by $\ell$.  The following lemma, whose proof we leave to the reader, then yields arcs $\ell_p$ and $\ell_q$ in $B_p$ and $B_q$, respectively, such that $\lambda_* \cup \ell_p \cup \ell_* \cup \ell_q$ is a smooth, simple closed curve that contains $\lambda$ and that is contained in $\Omega$.

\begin{lem}
Let $B$ be an open ball in $\R^n$, $n\geq 2$, let $a$ and $c$ be points outside $\overline B$, let $b$ and $d$ be points inside $B$.  Let $ab$ and $cd$ denote arcs of class $\sC^s$, $1\leq s\leq\infty$, the first with end points $a$ and $b$, the second with end points $c$ and $d$.  Suppose that $ab$ and $cd$ each intersect $bB$ transversally.  Let $a^*$ and $c^*$ be the unique points of $bB$ such that the subarcs $aa^*$ and $cc^*$ of $ab$ and $cd$, respectively, are disjoint from $B$.  Then there is an arc $a^*c^*$ with end points $a^*$ and $c^*$ contained in $B$ such that $aa^* \cup a^*c^*\cup cc^*$ is an arc of class $\sC^s$.
\end{lem}

{\bf Acknowledgment\/}.
This research was begun while the first author was a visitor at the University of Michigan.  He  thanks the Department of Mathematics for its hospitality.

\bibliography{Izzo.Stout.convergence}

\bigskip

\noindent Alexander Izzo, Department of Mathematics and Statistics, Bowling Green State University, Bowling Green, Ohio 43403

\noindent email: aizzo@bgsu.edu
\medskip

\noindent Edgar Lee Stout, Department of Mathematics, University of Washington, Seattle, Washington 98195

\noindent email: edgar.lee.stout@gmail.com

\bigskip\bigskip

\noindent \emph{2020 Mathematics Subject Classification\/}. Primary 32E20.

\noindent \emph{Key Words and Phrases\/}. polynomial convexity, polynomial hull, arc, curve, dense, open, rectifiable, smooth, Hausdorff metric.

\end{document}